\documentclass{amsart}
\usepackage{amssymb}
\usepackage{graphicx}
\usepackage[all]{xy}
\usepackage{a4wide}
\usepackage{txfonts}
\usepackage{hyperref}
\renewcommand{\theenumi}{\alph{enumi}}
\newtheorem{lemma}{Lemma}[section]
\newtheorem{theorem}[lemma]{Theorem}
\newtheorem{mainthm}{}
\renewcommand{\themainthm}{Main Theorem}
\newtheorem*{thm}{Theorem}
\newtheorem{corollary}[lemma]{Corollary}
\newtheorem{fact}[lemma]{Fact}
\newtheorem{proposition}[lemma]{Proposition}
\theoremstyle{remark}

\newtheorem{remark}[lemma]{Remark}

\newtheorem{example}[lemma]{Example}

\theoremstyle{definition}

\newtheorem{definition}[lemma]{Definition}
\newcounter{Lcounter}
\renewcommand{\theLcounter}{\theequation\alph{Lcounter}}

\newenvironment{labellist}{
    \refstepcounter{equation}
    \begin{list}{{\rm(\theLcounter)}}{\usecounter{Lcounter}\leftmargin=5pt}
    }{\end{list}
}
\def\freeprod{
\def\freeprodD{{\prod\kern-10.9pt{*}\kern5pt}}
\def\freeprodT{\prod\kern-9pt{*}\kern2pt}
\mathop{\mathchoice{\freeprodD}{\freeprodT}
{not defined}{not defined}}
}

\def\wre{\mathop{\rm wr}\nolimits}
\renewcommand{\wr}{\wre}

\def\Mhat{{\hat M}}
\def\Fhat{{\hat F}}

\newcommand{\calM}{{\mathcal M}}
\newcommand{\calN}{{\mathcal N}}
\newcommand{\calE}{{\mathcal E}}
\newcommand{\calX}{{\mathcal X}}

\newcommand{\Agag}{{\bar A}}
\newcommand{\bbZ}{\mathbb{Z}}
\newcommand{\bbN}{\mathbb{N}}
\newcommand{\bbF}{\mathbb{F}}
\newcommand{\bbP}{\mathbb{P}}
\newcommand{\bbQ}{\mathbb{Q}}
\def\alphagag{{\bar \alpha}}
\renewcommand{\phi}{\varphi}
\def\phigag{{\bar \varphi}}
\def\psigag{{\bar \psi}}

\def\st{\mathop{\mid\;}}
\newcommand{\isom}{\cong}
\newcommand{\gal}{{\rm Gal}}
\newcommand{\nek}{,\ldots,}
\def\normal{\triangleleft}
\newcommand{\Ker}{{\rm Ker}}
\newcommand{\rank}{{\rm rank}}
\newcommand{\weight}{{\rm weight}}
\newcommand{\Ind}{{\rm Ind}}
\newcommand{\cont}{\subseteq}
\newcommand{\hefresh}{\smallsetminus}
\newcommand{\pr}{{\rm pr}}
\newcommand{\id}{{\rm id}}
\newcommand{\Spec}{{\rm Spec}}
\newcommand{\ab}{{\rm ab}}

\begin{document}

\CompileMatrices 

\title[Semi-free profinite groups]{Permanence criteria for semi-free profinite groups}%
\author[Bary-Soroker]{Lior Bary-Soroker$^*$}
\address{School of Mathematical Sciences, Tel Aviv University, Ramat Aviv, Tel
Aviv, 69978, ISRAEL}%
\email{barylior@post.tau.ac.il}%
\thanks{$^*$ Supported in part by the Minkowski Center for Geometry at Tel Aviv University, established
by the Minerva Foundation and by the Israel Science Foundation grant 343/07.}

\author[Haran]{Dan Haran$^*$}
\address{School of Mathematical Sciences, Tel Aviv University, Ramat Aviv, Tel
Aviv, 69978, ISRAEL}%
\email{haran@post.tau.ac.il}%

\author[Harbater]{David Harbater$^{**}$}%
\address{Department of Mathematics, University of Pennsylvania, Philadelphia, PA
19104-6395, USA}%
\email{harbater@math.upenn.edu}%
\thanks{$^{**}$ Supported in part by NSF Grant DMS-0500118}

\date{\today}%
\dedicatory{Dedicated to Moshe Jarden on the occasion of his 65th birthday}%

\begin{abstract}
We introduce the condition of a profinite group being semi-free,
which is more general than being free and more restrictive than
being quasi-free.  In particular, every projective semi-free
profinite group is free.  We prove that the usual permanence
properties of free groups carry over to semi-free groups.  Using
this, we conclude that if $k$ is a separably closed field, then many
field extensions of $k((x,y))$ have free absolute Galois groups.
\end{abstract}

\maketitle

\section{Introduction and results}
A central problem is Galois theory is to understand the absolute
Galois groups of fields, and a key aspect is to find fields with
free absolute Galois groups. For example, if $C$ is an algebraically
closed field, then $K = C(x)$ is such a field.  This was proved for
$C = \mathbb C$ by Douady; and in the general case by Pop
\cite{Pop1995} and the third author \cite{Harbater1995}, with
another proof later by Jarden and the second author
\cite{HaranJarden2000}. The major conjecture in this context,
Shafarevich's conjecture, asserts that the maximal abelian extension
$\mathbb Q^{\ab}$ of the rational numbers $\mathbb Q$ has a free
absolute Galois group.

In \cite{HarbaterStevenson2005}, the third author and K.~Stevenson suggest a strategy for
proving the freeness of a profinite group: breaking the argument
into two simpler pieces, viz.\ quasi-freeness and projectivity.
This strategy was carried out in \cite{Harbater2006} in the context
of a two-variable Laurent series field $K = k((x,y))$.  For any base
field $k$, the absolute Galois group $\gal(K)$ is quasi-free
\cite{HarbaterStevenson2005}, though it is not free since it is not
projective.  In \cite{Harbater2006} the third author proves that the
commutator subgroup of a quasi-free group is quasi-free, and hence
$\gal(K^{\ab})$ is quasi-free. Now, if in addition $k$ is separably
closed, then $\gal(K^{\ab})$ is also projective. Therefore
$\gal(K^{\ab})$ is free, for such $k$. This can be viewed as an
analog of Shafarevich's conjecture.

In the above situation, it is key that the commutator subgroup of a
quasi-free group is quasi-free.  This leads to the question of when
a closed subgroup of a quasi-free group is quasi-free, particularly
in the case of projective subgroups.  Since closed subgroups inherit
projectivity, this question generalizes the corresponding classical
question about free subgroups of a free profinite group.  A partial
answer is given in \cite{RibesZalesskiiStevenson2006}, where Ribes,
Stevenson, and Zalesskii prove that an open subgroup of a quasi-free
group is quasi-free.

The classical question --- when is a closed subgroup of a free group
itself free --- has been dealt with in numerous papers, e.g.\
\cite{Haran1999JGroupTheory,Jarden2006GroupTheory,JardenLubotzky1992,JardenLubotzky1999,Melnikov1978}.
The second author 
has used twisted wreath products in \cite{Haran1999JGroupTheory} to
attack this question. Not only does this approach reprove many of
the previously known results, but it also proves the so-called
`Diamond Theorem' (see \cite[Theorem~25.4.3]{FriedJarden2005}):

\begin{thm}
Let $F$ be a free profinite group of infinite rank $m$. Let
$M_1,M_2$ be normal subgroups of $F$ and let $M$ be a subgroup of
$F$ such that $M_1\cap M_2\leq M$ but $M_1\not\leq  M$ and
$M_2\not\leq M$. Then $M$ is free of rank $m$.
\end{thm}
\noindent(The diagram
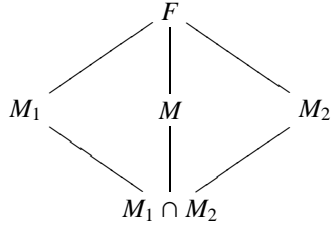
\[
\xymatrix{%
    &F\ar@{-}[dr]\ar@{-}[d]\\
M_1\ar@{-}[ur]\ar@{-}[dr]
    &M\ar@{-}[d]
        &M_2\\
    &M_1\cap M_2\ar@{-}[ur]
}%
\]
suggests the name Diamond Theorem.)
Recently the first author 
proved this theorem for finite $m\geq 2$ \cite{Bary-Soroker2006}.

It would thus be desirable to carry over this and other permanence
properties of free profinite groups to the class of quasi-free
profinite groups. However, our methods seem to work well only after
a slight modification of the notion: We say that a profinite group
of infinite rank $m$ is \textit{semi-free} if every nontrivial
finite split embedding problem for it has $m$ \textit{independent}
proper solutions.  (See Section~\ref{sec:def} below.)

The modified notion is in some ways more natural. First we have
\begin{enumerate}
\item
infinitely generated free profinite groups are semi-free (Theorem~\ref{thm QF and Proj imply
Free}),
\item
semi-free groups are quasi-free, but not vice-versa
(Proposition~\ref{prop QF not SQF}), \ and
\item
the absolute Galois group of $k((x,y))$ is semi-free
(Theorem~\ref{thm laurent series}).
\end{enumerate}

Moreover, we are able to prove the following theorem (where
case~\ref{diamond} corresponds to the Diamond Theorem above).  Also,
as Example~\ref{example qf no permanence} below shows, not all of
these properties hold for the class of quasi-free groups.

\begin{mainthm}\label{thm subgroups}
Let $F$ be a semi-free profinite group of infinite rank $m$ and let
$M$ be a closed subgroup of $F$. Then, in each of the following
cases the group $M$ is semi-free of rank $m$.
\begin{enumerate}
\renewcommand{\theenumi}{\Roman{enumi}}

\item \label{opensubgroup}
$(F:M) < \infty$.

\item \label{finitelygenerated}
$F/\Mhat$ is finitely generated, where $\Mhat = \bigcap_{\sigma \in
F} M^\sigma$ is the normal core of $M$.

\item\label{weight}
$\weight(F/M) < m$ (the definition of weight is recalled at Section~\ref{sec:weight}).

\item\label{weissauer}
$M$ is a proper subgroup of finite index of a closed normal subgroup
of $F$.

\item \label{commutator}
$M$ is normal in $F$, and $F/M$ is abelian.

\item \label{diamond}
There exist closed normal subgroups $M_1$, $M_2$ of $F$ such that
$M_1\cap M_2\leq M$ but $M_1\nleq M$ and $M_2\nleq M$.

\item \label{pronilpotent}
$M$ contains a closed normal subgroup $N$ of $F$ such that $F/N$ is
pronilpotent and $(F:M)$ is divisible by at least two primes.

\item\label{sparse}
$M$ is sparse in $F$ (see Definition~\ref{def_sparse}).

\item\label{finite exponent}
$(F:M) = \prod p^{\alpha(p)}$, where $\alpha(p)<\infty$ for all $p$.
\end{enumerate}
\end{mainthm}

The proof of \ref{thm subgroups} is in
Section~\ref{sec:proofofmainthm}.

This theorem gives rise to new constructions of fields having free
absolute Galois groups; see Section~\ref{sec free absolute}. One of
them generalizes the construction of fields with free absolute
Galois groups discussed above in the second paragraph of the
introduction. Another was provided by Jarden, using ideas of Pop.

We conclude the introduction with some ideas of the proof. The goal is to prove that $M$ is semi-free, i.e.\ that an arbitrary finite split embedding problem $\calE_1$ for $M$ has many independent proper solutions. We know that $M$ is a subgroup of a semi-free group $F$, so we wish to transfer the solvability problem to $F$. The first thing we do is to induce a split embedding problem $\calE$ for $F$ with the property that a weak solution of $\calE$ induces a weak solution to $\calE_1$ (see Proposition~\ref{thm main} for the exact definition of $\calE$). The embedding problem $\calE$ is constructed using a \emph{twisted wreath product} (see Definition~\ref{def twisted wreath product}). 

Now $\calE$ has many independent proper solutions because $F$ is semi-free. Each one of these proper solutions, say $\psi$, induces a solution $\nu$ of $\calE_1$. (Here $\nu = \pi\circ\psi|_{M}$, where $\pi$ is the Shapiro map; see Definition~\ref{def QF}.) We encounter two difficulties: (1) $\nu$ is not necessarily a \textit{proper} solution; (2) for two distinct proper solutions $\psi_1\neq \psi_2$ of $\calE$ we may get that $\nu_1= \nu_2$. 

We extract from \cite{Haran1999JGroupTheory} a condition under which $\nu$ remains a proper solution. This settles the first difficulty.  To treat (2), we use that fact that in our situation, $\psi_1, \psi_2$ are not only distinct, but also independent. Hence the image of $\psi_1\times \psi_2$ is also a wreath product (Lemma~\ref{lem twisted fiber}). This fact leads us to generalize the work in \cite{Haran1999JGroupTheory}, and find a necessary conditions for any two independent proper solutions $\psi_1,\psi_2$ to induce independent proper solutions $\nu_1,\nu_2$, as needed for $M$ to be semi-free.  See Proposition~\ref{thm main}\,\ref{induce independent}.

Note that this strategy does not apply to the corresponding problem for quasi-free groups, where the distinct proper solutions for a split embedding problem need not be independent, and since the image of $\psi_1\times \psi_2$ for distinct solutions $\psi_1,\psi_2$ of $\calE$ need not be a twisted wreath product in the absence of independence.  By avoiding this difficulty, our focus on semi-free groups permits us to show that many subgroups of semi-free groups are semi-free (and in particular quasi-free); and that if such a subgroup is also projective then it is free (see Theorem~\ref{thm QF and Proj imply Free}).

\section{Independent subgroups and solutions of embedding problems} \label{sec:def}
\begin{definition}
Let $F$ be a profinite group.
\begin{enumerate}
\item
Open subgroups $M_1,\ldots,M_n$ of $F$ are \textbf{$F$-independent}
if
$$
(F : \bigcap_{i=1}^n M_i) = \prod_{i=1}^n (F : M_i) .
$$
If $M_1,\ldots,M_n$ are normal in $F$, this is equivalent to
$$F/\bigcap_{i=1}^n M_i \isom \prod_{i=1}^n F/M_i$$
\item
A family $\calM$ of open subgroups of $F$ is
\textbf{$F$-independent} if every finite subset of $\calM$ is $F$-independent.
\end{enumerate}
\end{definition}

The notion of $F$-independence coincides with independence with
respect to the Haar probability measure on $F$
\cite[Section~18.3]{FriedJarden2005}.  There is also the following equivalent characterization of independence:  Open subgroups $M_1, \ldots, M_n$ are $F$-independent if and only if $F$ acts transitively on $\prod_{i=1}^n F/M_i$.  This criterion can be used to obtain alternative short proofs of parts~\ref{inductive independent} and \ref{intermediate} in Proposition~\ref{proposition independent} below.

A key example of independence occurs in the case of a Galois field extension $L/K$.  If $F = \gal(L/K)$ and $L_1,\ldots,
L_n$ are the fixed fields of $M_1,\ldots,M_n$ in $L$, then by the
Galois correspondence, $M_1,\ldots,M_n$ are $F$-independent if and
only if $L_1,\dots, L_n$ are linearly disjoint over $K$.

The following properties can be either proven directly or deduced
from the corresponding properties of linear disjointness of fields:

\begin{proposition}\label{proposition independent}
Let $M_1\nek M_n$ be open subgroups of a profinite group $F$.
\begin{enumerate}

\item
$(F : \bigcap_{i=1}^n M_i) \le \prod_{i=1}^n (F : M_i)$.

\item
Let $M_1 \le N_1 \le F$. Then $M_1,M_2$ are $F$-independent if and
only if $N_1,M_2$ are $F$-independent and $M_1,N_1 \cap M_2$ are
$N_1$-independent.

\item \label{inductive independent}
The subgroups $M_1,\ldots,M_n$ are $F$-independent if and only if
$M_1,\ldots,M_{n-1}$ are $F$-independent and $\bigcap_{i=1}^{n-1}
M_i$, $M_n$ are $F$-independent.

\item \label{intermediate}
Let $M_i \le N_i \le F$ for each $1 \le i \le n$. If
$M_1,\ldots,M_n$ are $F$-independent, then so are $N_1,\ldots,N_n$.

\item
Suppose $M_1 \normal F$. Then $M_1,M_2$ are $F$-independent if and
only if $F = M_1 M_2$.
\end{enumerate}
\end{proposition}

\begin{proof}
(a) This follows by induction from the case $n=2$, which is standard.

(b) First assume $M_1,M_2$ are $F$-independent. Then, since $(N_1 \cap M_2 : M_1\cap M_2) \leq (N_1 : M_1)$ we have
\begin{eqnarray*}
 (F: N_1 \cap M_2) &=& \frac{(F:M_1\cap M_2)}{(N_1\cap M_2 : M_1\cap M_2)} =  \frac{(F:M_1)(F:M_2)}{(N_1\cap M_2 : M_1\cap M_2)}\\
 	& =&\frac{(F:N_1)(N_1:M_1)(F:M_2)}{(N_1\cap M_2 : M_1\cap M_2)}\geq 
		(F:N_1)(F:M_2).
\end{eqnarray*}
Therefore equality holds by (a), and $N_1,M_2$ are $F$-independent. Similarly, since $(N_1:N_1\cap M_2) \leq (F:M_2)$ we have 
\begin{eqnarray*}
(N_1 : M_1\cap (N_1 \cap M_2))  
	&= 
		&\frac{(F:M_1\cap M_2)}{(F:N_1)}=\frac{(F:M_1)(F:M_2)}{(F:N_1)}\\
	&\geq 
		&(N_1:M_1) (N_1:N_1\cap M_2),
\end{eqnarray*}
so $M_1, N_1\cap M_2$ are $N_1$-independent by (a).  Conversely, 
\begin{eqnarray*}
(F: M_1\cap M_2) 
	&=
		& (F:N_1) (N_1 : M_1 \cap (N_1 \cap M_2)) = (F : M_1) (N_1 : N_1\cap M_2) \\
	&=
		& (F:M_1) \frac{(F:N_1\cap M_2)}{(F:N_1)}=(F:M_1) (F:M_2).
\end{eqnarray*}

(c) By part (a), 
\[(F : \bigcap_{i=1}^n M_i) \le (F : \bigcap_{i=1}^{n-1} M_i) (F:M_n) \le \prod_{i=1}^n (F : M_i).\] 
So $(F : \bigcap_{i=1}^n M_i) = \prod_{i=1}^n (F : M_i)$ if and only if the above two inequalities are equalities, and the assertion follows.

(d) Since $(\bigcap_i M_i : \bigcap_i N_i) \leq  \prod_i (M_i:N_i)$ we have 
\[
(F: \bigcap_i N_i) =\frac{(F:\bigcap_i M_i)}{(\bigcap_i M_i : \bigcap_i N_i)} \geq \frac{\prod_i(F:M_i)}{\prod_i(M_i:N_i)} = \prod_i (F:N_i),
\]
so equality holds by (a). 

(e) We have $(M_1M_2 : M_1) = (M_2 : M_1 \cap M_2)$. Thus
\begin{eqnarray*}
(F : M_1) (F : M_2) &=& (F : M_1M_2) (M_2 : M_1 \cap M_2) (F :
M_2)\\
&=& (F : M_1M_2) (F : M_1 \cap M_2) .
\end{eqnarray*}
\end{proof}

Recall that an \textbf{embedding problem} for a profinite group $F$
is a pair of epimorphisms of profinite groups
\begin{equation}
\label{EP for F}%
(\phi\colon F\to G, \alpha\colon H\to G).
\end{equation}
The embedding problem is called \textbf{finite} if $H$ and $G$ are
finite. It is called \textbf{split} (respectively
\textbf{nontrivial}) if $\alpha$ splits (respectively is not an
isomorphism). We abbreviate `finite split embedding problem' and
write `FSEP'. A \textbf{(weak) solution} for an embedding problem is
a homomorphism $\psi\colon F\to H$ with $\alpha\circ\psi=\phi$. A
solution is said to be \textbf{proper} if it is surjective.

\begin{definition}
We call solutions of a finite embedding problem \eqref{EP for F}
\textbf{independent} if their kernels are $\Ker \phi$-independent.
\end{definition}

We now introduce a criterion for the independence of proper
solutions of finite embedding problems in terms of fiber products of
groups.

Let $\{\alpha_i\colon H_i\to G \st i \in I\}$ be a family of
epimorphisms of profinite groups. Their \textbf{fiber product} with
respect to the $\alpha_i$'s is defined by
\[
{\varprod}_{G}H_i = \Big\{h \in \prod H_i \st \alpha_i(h_i) = \alpha_j(h_j)
\ \forall i,j \in I \Big\}.
\]
(Here $h_i = h(i)$ is the value of $h$ at $i$.) This is a closed
subgroup of $\prod H_i$, hence a profinite group. The projection on
the $i$-th coordinate, $\pr_i\colon \varprod_G H_i \to H_i$, is
surjective. The fiber product is equipped with a canonical
epimorphism $\alpha^I = \alpha_i\circ \pr_i \colon \varprod_G H_i \to
G$, which is independent of $i \in I$.

In particular, if $I$ is a finite set, say $I = \{1,\ldots, n\}$,
then
\[
{\varprod}_G H_i = H_1\times_G \cdots \times_G H_n = \{(h_1,\cdots,
h_n) \in \prod H_i \st \alpha_1(h_1) = \cdots = \alpha_n(h_n)\}.
\]

Fiber products are associative:

\begin{lemma}\label{lem associative}
Let $\alpha_i\colon H_i \to G_0$, $i=1\nek n$, and $\beta\colon G\to
G_0$ be epimorphisms of finite groups. Then the natural map $\big(\varprod_{G_0} H_i\big)
\times_{G_0} G \to \varprod_{G} (H_i\times_{G_0} G)$ is an isomorphism.
\end{lemma}

\begin{proof} An element in $\big(\varprod_{G_0} H_i\big)\times_{G_0} G$
is of the form $((h_1, \ldots, h_n), g)$, where the elements $h_i \in H_i$ and $g \in G$ all have the same image in $G_0$.  An element in $\varprod_{G} (H_i\times_{G_0} G)$ is of the form $((h_1,g) \ldots, (h_n, g))$, for such elements $h_i \in H_i$ and $g \in G$, because the fiber product is taken over $G$.  The
map that takes $((h_1, \ldots, h_n), g)$ to $((h_1,g) \ldots, (h_n, g))$ is clearly an isomorphism. 
\end{proof}

A key property, in our setting, of fiber products is that solutions
$\psi_i$ of embedding problems $(\phi\colon F\to G, \alpha_i \colon
H_i\to G)$, $i \in I$, induce a canonical solution, $\psi^I = \prod
\psi_i$, of the embedding problem $(\phi\colon F\to G, \alpha^I
\colon \varprod_{G} H_i\to G)$. More precisely, $(\psi^I(x))_i =
\psi_i(x)$ for each $x \in F$; e.g., if $I = \{1,\ldots, n\}$, then
$\psi^I(x) = (\psi_1(x),\cdots,\psi_n(x))$. We obtain the original
solutions via the projection on the coordinates, i.e. $\psi_i =
\pr_i \circ \psi^I$ for each $i \in I$. In particular, taking $F =
G$ and $\phi = \id$, we see that if all the $\alpha_i$'s split, so
does $\alpha^I$.

Given a single epimorphism $\alpha\colon H\to G$ and a set $I$, we
write $H_G^I$ for the fiber product $\varprod_G H_i$, where $H_i = H$
and $\alpha_i = \alpha$ for each $i \in I$.

\begin{lemma}
\label{lem_fiber independent}%
Let $I$ be a set and let $\calE = (\phi\colon F \to G, \alpha\colon
H\to G)$ be a finite embedding problem for a profinite group $F$.
Put $\calE^I = (\phi\colon F\to G, \alpha^I \colon H^I_G\to G)$.
Then solutions $\{\psi_i\}_{i \in I}$ of $\calE$ are independent and
proper if and only if the solution $\psi^I = \prod \psi_i$ of
$\calE^I$ is proper.
\end{lemma}

\begin{proof}
We first assume that $I$ is finite, $I = \{1,\ldots, n\}$. If one of
the $\psi_i$'s is not surjective, then $\psi^I$ is not surjective.
Hence, we may assume that $\psi_1\nek \psi_n$ are surjective. Let $K
= \Ker\phi$ and $M_i = \Ker \psi_i$, $i=1\nek n$. By the definition
of $\psi^I$ we have $\Ker\psi^I = \bigcap_{i=1}^n M_i$. Since
$|H^I_G| = |H|^n/|G|^{n-1}$, we get that $\psi^I$ is surjective if
and only if $(F:\bigcap_{i=1}^n M_i) = |H|^n/|G|^{n-1}$. But
$(F:\bigcap_{i=1}^n M_i) = (F:K)(K:\bigcap_{i=1}^n
M_i)=|G|(K:\bigcap_{i=1}^n M_i)$; hence $\psi^I$ is surjective if
and only if $(K:\bigcap_{i=1}^n M_i) = |H|^n/|G|^n = \prod_{i=1}^n
(K:M_i)$, as desired.

In the general case $H^I_G$ is the inverse limit of $H^J_G$, where
$J$ runs through the finite subsets of $I$ and the epimorphisms
$\pr^J \colon H^I_G \to H^J_G$ are given by the restriction of
coordinates from $I$ to $J$. Obviously, $\psi^J = \pr^J \circ
\psi^I$, for each $J$. Hence $\psi^I$ is proper if and only if all
$\psi^J$'s are proper. By the first paragraph of this proof this
happens if and only if the $\psi_i$'s are independent and proper.
\end{proof}

\section{Semi-free profinite groups}

\begin{definition}
A profinite group $F$ of infinite rank is \textbf{quasi-free} if there exists an
infinite cardinal $m$ such that every nontrivial FSEP for $F$ has
exactly $m$ distinct proper solutions (see
\cite{Harbater2006,HarbaterStevenson2005,RibesZalesskiiStevenson2006}).
By~\cite[Lemma~1.2]{RibesZalesskiiStevenson2006} such a group is necessarily of rank~$m$.
\end{definition}

In the following definition we give a stronger variant of
quasi-freeness.

\begin{definition}
\label{def QF} A profinite group $F$ of infinite rank is
\textbf{semi-free}\footnote{a term coined by Moshe Jarden as an
alternative to ``strongly quasi-free'', which we initially used.} if every nontrivial FSEP for $F$ has $m$ independent proper solutions, where $m$ is the rank of $F$.
\end{definition}

\begin{remark} The above definitions consider only \textit{infinitely} generated profinite groups, with the notions of quasi-free and semi-free being left undefined in the finitely generated case.  The reason is that for a profinite group $F$ of finite rank $m$, there is no proper solution to \textit{any} finite embedding problem $\calE = (\phi\colon F \to G, \alpha\colon H\to G)$ for which $H$ has rank greater than $m$.  By leaving the notions undefined in the finitely generated case, we thus avoid the perverse situation in which a finitely generated free group would violate the conditions of being quasi-free or semi-free.  
One could instead consider the class of groups $F$ of finite rank for which there is a proper solution to every FSEP $\calE$ for which $\rank(H) \le \rank(F)$.  But a finite rank group would satisfy that condition if and only if it is free, by \cite[Lemma~17.7.1]{FriedJarden2005}; so this would not be a new condition on such groups.  For the purposes of this paper, the case of infinite rank is sufficient to consider, and we restrict to that situation. 
\end{remark}

\begin{remark}
In Definition~\ref{def QF}, it would suffice to assume just that rank $F$ is at most $m$.  More precisely, let $F$ be a profinite group and let $m$ be an infinite cardinal.
Assume that $\rank\ F \leq m$ and every nontrivial FSEP for $F$ has
$m$ independent proper solutions. Then $\rank\ F = m$, and thus $F$
is semi-free.

Indeed, consider any nontrivial FSEP and let $\{\psi_i \mid i<m\}$ be a
set of independent proper solutions. Then $\Ker \psi_i \neq \Ker
\psi_j$ for all $i\neq j$. This implies that $F$ has at least $m$
open subgroups, the set $\{\Ker \psi_i \mid i<m\}$, and hence $\rank\
F\geq m$ (see \cite[Proposition 17.1.2]{FriedJarden2005}). Therefore
$\rank\ F = m$, as needed.
\end{remark}

Clearly, every semi-free group is quasi-free. One might suspect that
the opposite is also true. If $m=\aleph_0$, then for both notions it
suffices to have one proper solution of any nontrivial FSEP (see the
lemma below), and hence they are equivalent. If $m>\aleph_0$, then
there are quasi-free groups that are not semi-free. We postpone the
discussion of this to Section~\ref{sec:qfvssqf}.

\begin{lemma}\label{lem:rankaleph0}
Let $F$ be a countably generated profinite group. Then $F$ is
semi-free of rank $\aleph_0$ if and only if every FSEP for $F$ is properly
solvable.
\end{lemma}

\begin{proof}
Let $\calE = (\phi_0\colon F\to G, \alpha_0\colon H\to G)$ be a
nontrivial FSEP. For each integer $n>0$, let $\alpha_{n-1} \colon
H^n_G \to H^{n-1}_G$ be the projection map. Inductively, we can find
solutions $\phi_n \colon F\to H^n_G$ of the FSEP
\[
\calE_n = (\phi_{n-1} \colon G \to H^{n-1}_G, \alpha_{n-1}\colon
H^n_G\to H^{n-1}_G).
\]
Then $\phi := \varprojlim\phi_n \colon G \to
H^{\bbN}_G$ is surjective. Lemma~\ref{lem_fiber independent}
implies the existence of $\aleph_0$ independent proper solutions, and thus
$F$ is semi-free.
\end{proof}

We extend \cite[Theorem 2.1]{HarbaterStevenson2005}:

\begin{theorem}\label{thm QF and Proj imply Free}
Let $F$ be a profinite group of infinite rank $m$. The following
conditions are equivalent:
\begin{enumerate}
\item \label{QFProj a}
$F$ is free.
\item \label{QFProj b}
$F$ is semi-free and projective.
\item \label{QFProj c}
$F$ is quasi-free and projective.
\end{enumerate}
\end{theorem}

\begin{proof}
We show that \eqref{QFProj a} $\Rightarrow$ \eqref{QFProj b}. Let
$\calE = (\phi\colon F\to G,\alpha\colon H \to G)$ be a nontrivial
finite embedding problem for $F$. Fix a set $I$ of cardinality $m$.
Let $H^I_G$ be the corresponding fiber product; let $\pr_i\colon
H^I_G \to H$ be the projection on the $i$-th coordinate, for each $i
\in I$; and let $\alpha^I = \alpha \circ \pr_i \colon H^I_G \to G$
be the canonical epimorphism.

Since $F$ is free of rank $m$ and since $\rank(H^I_G) \leq m$, we
have a proper solution $\psi\colon F \to H^I_G$ of the embedding
problem $(\phi\colon F\to G, \alphagag\colon H^I_G \to G)$
\cite[Theorem~3.5.9]{RibesZalesskii2000}. Put $\psi_i =
\pr_i\circ\psi$ for each $i \in I$. Then, by Lemma~\ref{lem_fiber
independent}, the solutions $\{\psi_i\}_{i \in I}$ of $\calE$ are
independent and proper. As $\calE$ is nontrivial, they are distinct.

Implication \eqref{QFProj b} $\Rightarrow$ \eqref{QFProj c} is
trivial and \eqref{QFProj c} $\Rightarrow$ \eqref{QFProj a} is
\cite[Theorem 2.1]{HarbaterStevenson2005}.
\end{proof}

From a technical point of view, it is preferable to work with a set of
\emph{pairwise} proper solutions of a FSEP instead of independent
set of solutions. The following result shows that it is possible.

\begin{proposition}
Let $\calM$ be an infinite family of pairwise $F$-independent open
normal subgroups of a profinite group $F$. Then $\calM$ contains an
$F$-independent subfamily $\calM_0$ of cardinality $|\calM|$.
\end{proposition}

\begin{proof}
By Zorn's Lemma there is a maximal $F$-independent subfamily
$\calM_0$ of $\calM$. We have to show that $|\calM_0| = |\calM|$.
Assume the contrary; that is, $|\calM_0| < |\calM|$.

Let $\calM_1$ be the family of all finite intersections of the
elements of $\calM_0$. If $\calM_0$ is finite, then so is $\calM_1$;
if $\calM_0$ is infinite, then $|\calM_1| = |\calM_0|$. In
particular, $|\calM_1| < |\calM|$. The groups in $\calM_1$ are open
in $F$. Let $\calM_2$ be the family of all open subgroups of $F$
containing a group in $\calM_1$. Again, if $\calM_1$ is finite, then
so is $\calM_2$; if $\calM_1$ is infinite, then $|\calM_2| =
|\calM_1|$. In particular, $|\calM_2| < |\calM|$.

For every proper subgroup $N$ of $F$ there exists at most one $M \in
\calM$ such that $M \le N$. Indeed, if $M_1, M_2 \in \calM$ are
distinct, then $M_1M_2 = F$, by Proposition~\ref{proposition
independent}(e), and hence we cannot have $M_1,M_2 \le N < F$. Since
$|\calM_2| < |\calM|$, there exists $M \in \calM$ such that
{\renewcommand{\theequation}{*}\addtocounter{equation}{-1}
\begin{eqnarray}\label{eq condition *}
M \le N \in \calM_2 \ \ \mbox{only for}\ \ N=F.
\end{eqnarray}}%
We claim that $\calM_0\cup\{M\}$ is $F$-independent. (This will
produce the desired contradiction to the maximality of $\calM_0$.)
Thus we have to show, for distinct $M_1, \ldots, M_n \in \calM_0$,
that $M_1, \ldots, M_n, M$ are $F$-independent.

Put $N = \bigcap_{i=1}^n M_i$. By Proposition~\ref{proposition
independent}(c) it suffices to show that $M,N$ are $F$-independent.
By construction, $N \in \calM_1$. Hence $MN \in \calM_2$. Since $M
\le MN$, by \eqref{eq condition *}, $MN = F$. Hence, by
Proposition~\ref{proposition independent}(e), $M,N$ are
$F$-independent.
\end{proof}

\begin{corollary}
\label{cor pairwise}%
Let $m$ be an infinite cardinal and let $F$ be a profinite group of
rank at most $m$. Then $F$ is semi-free of rank $m$ if and only if
every nontrivial FSEP has $m$ pairwise independent proper solutions.
\end{corollary}

\section{Finite split embedding problems and twisted wreath products}
We follow \cite{Haran1999JGroupTheory} and establish the connection
between FSEPs and twisted wreath products.

\begin{definition}
[Twisted wreath product] \label{def twisted wreath product} Let $A$,
$G_0\leq G$ be finite groups with a (right) action of $G_0$ on $A$.
Write $\Ind_{G_0}^G(A)$ for all functions $f\colon G \to A$ such
that $f(\sigma \tau) = f(\sigma)^\tau$ for all $\sigma\in G$ and
$\tau\in G_0$ with componentwise multiplication. Then $\Ind_{G_0}^G
(A) \cong A^{(G:G_0)}$ and $G$ acts on $\Ind_{G_0}^G(A)$ by
\[
f^\sigma(\rho) = f(\sigma\rho), \qquad \sigma,\rho\in G,
f\in\Ind_{G_0}^G(A).
\]
The \textbf{twisted wreath product}, $A\wr_{G_0} G$, is defined to
be  the semidirect product of $\Ind_{G_0}^G(A)$ and $G$, i.e.
$A\wr_{G_0} G = \Ind_{G_0}^G(A) \rtimes G$. Here and below, $\alpha:
A\wr_{G_0} G\to G$ denotes the canonical projection
$f\sigma\mapsto\sigma$ (see \cite[Definition
13.7.1]{FriedJarden2005}). Similarly, $\alpha_0: A \rtimes G_0 \to
G_0$ denotes the canonical projection $a\sigma \mapsto \sigma$ of
the semidirect product.

There is an epimorphism $\pi_0\colon \Ind_{G_0}^G(A) \to A$ defined
by $\pi_0(f) = f(1)$. It extends to an epimorphism $\pi\colon
\Ind_{G_0}^G(A)\rtimes G_0 \to A\rtimes G_0$ defined by
$f\tau\mapsto f(1)\tau$ for $f\in \Ind_{G_0}^G(A)$ and $\tau\in
G_0$, since $\pi_0(f^\tau)=f^{\tau}(1) =
f(\tau)=f(1)^\tau=\pi_0(f)^\tau$ for all $f\in \Ind_{G_0}^G(A)$ and
$\tau\in G_0$. We call $\pi$ the \textbf{Shapiro map} of $A\wr_{G_0}
G$.
\end{definition}

\begin{remark}
\begin{enumerate}
\item If  $G=G_0$ in Definition~\ref{def twisted wreath product}, then $A
\wr_{G_0} G = A \rtimes G$.
\item 
See \cite{RibesSteinberg}, where a related notion, known as a permutational wreath product, is used in a similar context.
\end{enumerate}
\end{remark}

The following technical result will be needed later.

\begin{lemma}
\label{lem large subgroup}%
Under the above notation, let $B = \pi^{-1}(G_0)$. Then $B$ is a
subgroup of $A\wr_{G_0} G$ of index $(G:G_0)|A|$. If $A\neq 1$, then
$B$ does not contain $\Ind_{G_0}^G(A)$.
\end{lemma}

\begin{proof}
As the Shapiro map $\pi$ is surjective, $(\Ind_{G_0}^G(A)\rtimes G_0
: B) =|A|$. Thus the index of $B$ in $A\wr_{G_0}G$ is $(G:G_0)|A|$.

If $A \ne 1$, there is $f \in \Ind_{G_0}^G(A)$ such that $f(1) \ne
1$; then $\pi(f) \notin G_0$, and hence $f \notin B$.
\end{proof}

\begin{lemma}
\label{lem twisted fiber}%
Consider groups $H_i = A_i\wr_{G_0} G$, for $i=1\nek n$. Then $G_0$
acts on $\prod A_i$ componentwise and $\varprod_{G}H_i \cong (\prod
A_i) \wr_{G_0} G$.
\end{lemma}

\begin{proof}
We have
\begin{eqnarray*}
{\varprod}_{G}H_i &=& \{\big( (f_1 \sigma),\ldots,(f_n \sigma) \big)
\st
f_i \in \Ind_{G_0}^G(A_i), \ \sigma \in G \} , \\
(\prod A_i) \wr_{G_0} G &=& \{(f_1,\ldots,f_n)\sigma \st f_i \in
\Ind_{G_0}^G(A_i), \ \sigma \in G \},
\end{eqnarray*}
and the isomorphism is given by $\big( (f_1 \sigma),\ldots,(f_n
\sigma) \big) \mapsto (f_1,\ldots,f_n)\sigma$.
\end{proof}

\begin{lemma}
\label{lem domination by wreath product}%
Let $\phi \colon F \to G$ be an epimorphism of a profinite group $F$
onto a finite group $G$. Let $M$ be a closed subgroup of $F$, let
$G_0 = \phi(M) \le G$, and assume that $G_0$ acts on a finite group
$A$. Consider the FSEP
$$
\calE_0(A) = (\phi|_{M}\colon M\to G_0, \alpha_0\colon A\rtimes
G_0\to G_0),
$$
and let $\psi$ be a solution of the corresponding FSEP
$$
\calE(A) = (\phi\colon F\to G, \alpha\colon A\wr_{G_0} G \to G),
$$
with notation as in Definition~\ref{def twisted wreath product}.
Let $\pi$ be the Shapiro map of $A\wr_{G_0} G$. Then $\psi(M) \leq
\Ind_{G_0}^G(A)\rtimes G_0$ and $\pi \circ \psi|_{M}$ is a solution
of $\calE_0(A)$.
\end{lemma}

\begin{proof}
We have $\psi(M) \le \alpha^{-1}(G_0) = \Ind_{G_0}^G(A)\rtimes G_0$.
Thus $\pi\circ \psi|_{M}$ is defined. Let $\alpha' \colon
\Ind_{G_0}^G(A)\rtimes G_0 \to G_0$ be the restriction of $\alpha$.
From the commutativity of
\[
\xymatrix@=10pt{%
        &&M\ar[dll]_{\psi|_M}\ar[d]^{\phi|_M}\\
\Ind_{G_0}^G(A) \rtimes G_0\ar[rr]^-{\alpha'}\ar[dr]_{\pi}
        &&G_0\\
    &A\rtimes G_0 \ar[ur]_{\alpha_0}
}%
\]
we have $\alpha_0 \circ \pi\circ \psi|_{M} = \phi|_{M}$, i.e.\
$\pi\circ \psi|_M$ is a solution.
\end{proof}

Although the solution $\pi \circ \psi|_{M}$ in the preceding lemma
need not be proper, even if $\psi$ is proper, the proof of
\cite[Proposition 25.4.1]{FriedJarden2005} shows that, under some
assumptions on $M$, the properness of $\psi$ does imply
the properness of $\pi \circ \psi|_{M}$. Moreover, if $F$ is a free
profinite group of infinite rank $m$, that proof produces a family
of $m$ distinct proper solutions of $\calE_0(A)$. We generalize this
in part~\ref{induce independent} of the following proposition, where we consider proper solutions that are not just distinct, but in fact independent.

\begin{proposition}
\label{thm main}%
Let $M\leq F$ be profinite groups, let $A,G_1$ be finite groups
together with an action of $G_1$ on $A$, and let
\[
\calE_1(A) =( \mu\colon M \to G_1, \alpha_1\colon A\rtimes G_1 \to
G_1)
\]
be a FSEP for $M$. Let $D, F_0, L$ be subgroups of $F$ such that
\begin{labellist}\label{eq mainthm}
\item \label{eq main thma}
$D$ is an open normal subgroup of $F$ with $M \cap D \le \Ker\mu$,
\item \label{eq main thmb}
$F_0$ is an open subgroup of $F$ with $M\leq F_0\leq MD$,
\item \label{eq main thmc}
$L$ is an open normal subgroup of $F$ with $L \leq F_0 \cap D$.
\end{labellist}
Put $G = F/L$, $G_0 = F_0/L \le G$, and let $\phi\colon F\to G$ be
the quotient map.
\begin{enumerate}
\item
\label{epi exists} Then there is an epimorphism $\phigag_1 \colon
G_0 \to G_1$, through which an action of $G_0$ on $A$ is defined,
such that every weak solution $\psi$ of the FSEP
\[
\calE(A) = (\phi\colon F\to G, \alpha\colon A\wr_{G_0} G \to G)
\]
induces a weak solution $\nu = \rho\circ \pi \circ \psi|_{M}$ of
$\calE_1(A)$. Here $\pi$ is the Shapiro map of $A\wr_{G_0} G$ and
$\rho \colon A\rtimes G_0 \to A\rtimes G_1$ is the extension of
$\phigag_1$ by the identity of $A$.
\item
\label{induce independent} Let $n \in \bbN$. Assume that there is a
closed normal subgroup $N$ of $F$ with $N \leq M \cap L$ such that
there is no nontrivial quotient $\Agag$ of $A^n$ through which the
action of $G_0$ on $A^n$ descends and for which the FSEP
\begin{equation}\label{eq twisted wreath product EP}
  (\phigag\colon F/N\to G, \alphagag\colon\Agag\wr_{G_0}G \to G),
\end{equation}
where $\phigag$ is the quotient map, is properly solvable.
Then any $n$ independent proper solutions $\psi$ of $\calE(A)$
induce $n$ independent proper solutions $\nu$ of $\calE_1(A)$.
\[
\xymatrix@=15pt{%
    &M\ar@{-}[r]
        &F_0\ar@{-}[r]
            &MD\ar@{-}[r]
                &F\\
    &M\cap D\ar@{-}[u]|{\ker \mu}\ar@{-}[r]
        &F_0\cap D\ar@{-}[r]\ar@{-}[u]
            &D\ar@{-}[u]\\
N\ar@{-}[r]
    &M\cap L\ar@{-}[u]\ar@{-}[r]
        &L\ar@{-}[u]
}%
\]
\end{enumerate}
\end{proposition}

\begin{proof}
(\ref{epi exists}) We can extend $\mu$ to a map $MD \to G_1$ by $md
\mapsto \mu(m)$ for all $m\in M$ and $d \in D$. Its restriction to
$F_0$ is an epimorphism $\phi_1 \colon F_0 \to G_1$. It decomposes
as $\phi_1 = \phigag_1\circ \phi_0$, where $\phi_0 \colon F_0 \to
G_0$ is the restriction of $\phi$ to $F_0$ and $\phigag_1 \colon G_0
\to G_1$ is an epimorphism. (Here we use that $\Ker \phi|_{F_0} = L
\le D \le \Ker \phi_1$ to obtain $\phigag_1$.) Let $G_0$ act on $A$
via $\phigag_1$. Then we have the following commutative diagram
\[
\xymatrix{%
    &F_0\ar[d]_{\phi_0}\ar@/^/[dd]^(.6){\phi_1} \ar[r] &F \ar[d]^\phi \\
A\rtimes G_0\ar[r]^(0.6){\alpha_0} \ar[d]_{\rho}
    &G_0\ar[d]_{\phigag_1} \ar[r] & G\\
A\rtimes G_1\ar[r]^(0.6){\alpha_1}
    &G_1,
}%
\]
where $\rho$ is given by $\rho|_{G_0} = \phigag_1$ and $\rho|_{A} =
\id_A$. By Lemma~\ref{lem domination by wreath product}, $\pi \circ
\psi|_{M}$ is a (not necessarily proper) solution of $\calE_0(A):
(\phi_0|_{M} \colon M \to G_0, \alpha_0 \colon A\rtimes G_0 \to
G_0)$. Hence $\nu = \rho\circ \pi \circ \psi|_{M}$ is a solution of
$\calE_1(A)$.

(\ref{induce independent}) Let $\{\psi_i\}_{i =1}^n$ be a family of
independent proper solutions of $\calE(A)$. Let $1\leq i\leq n$, and
let $\nu_i = \rho\circ \pi \circ \psi_i|_{M}$ be the induced
solution of $\calE_1(A)$, as in (\ref{epi exists}). It suffices to
show that each $\nu_i$ is proper and the family $\{\nu_i\}_{i =
1}^n$ is independent.

By Lemma~\ref{lem twisted fiber}, $(A\wr_{G_0} G)^n_G = A^n
\wr_{G_0} G$.  So by Lemma~\ref{lem_fiber independent}, $\psi_1\nek
\psi_n$ define a proper solution, $\psi\colon F\to A^n\wr_{G_0} G$,
of
\[
\calE(A^n) = (\phi\colon F\to G, \alpha \colon A^n\wr_{G_0} G \to
G).
\]
Applying Lemma~\ref{lem domination by wreath product}, with $A^n$
playing the role of $A$ there, we get that $\nu = \rho'\circ
\pi'\circ \psi$ is a solution of
\[
\calE_1(A^n) = (\mu\colon M\to G_1, \alpha_1 \colon A^n\rtimes
G_1\to G_1).
\]
(Here $\rho'$ and $\pi'$ are defined as $\rho$ and $\pi$ with $A^n$
replacing $A$.) By Part C of \cite[Proposition
25.4.1]{FriedJarden2005} (again, with $A^n$ replacing $A$),
$\pi'(\psi(N))=A^n$. But $\nu(N) = \rho'(\pi'(\psi(N))) =
\rho'(A^n)=A^n$. Therefore $A^n\leq \nu(M)$, and thus $\nu$ is a
proper solution of $\calE_1(A^n)$. As $\psi = \prod \psi_i$, we get
that $\nu = \prod \nu_i$. Consequently, $\nu_1\nek \nu_n$ are
independent proper solutions (Lemma~\ref{lem_fiber independent}).
\end{proof}

\begin{corollary}[cf. {\cite[Proposition 25.4.1]{FriedJarden2005}}]
\label{cor_main}%
Let $F$ be a semi-free profinite group of infinite rank $m$ and let
$M$ be a closed subgroup of $F$. Assume that for every open normal
subgroup $D$ of $F$ there exist $L$ and $F_0$ as in \eqref{eq main
thmb},\eqref{eq main thmc} of Proposition~\ref{thm main}, and there
exists $N \normal F$ with $N \leq M \cap L$ such that no FSEP
\[
(\phi \colon F/N \to F/L, \alpha\colon A\wr_{F_0/L} F/L \to F/L),
\]
where $A$ is a nontrivial finite group on which $F_0/L$ acts and
where $\phi$ is the quotient map, is properly solvable.

Then $M$ is semi-free of rank $m$.
\end{corollary}

\begin{proof}
By \cite[Corollary 17.1.4]{FriedJarden2005}, $\rank(M)\leq \rank(F)
= m$. Let $\calE_1(A)$ be a FSEP as in Proposition~\ref{thm main}.
Choose $D$ as in \eqref{eq main thma} of Proposition~\ref{thm main}. With $F_0,L, N$ be as above,
let $\calE(A)$ be as in Proposition~\ref{thm main}. Since $F$ is
quasi-free of rank $m$, there exists a family $\Psi$ of independent
proper solution of $\calE(A)$ of cardinality $m$. This in turn
induces a family $\calN$ of solutions of $\calE_1(A)$
(Lemma~\ref{lem domination by wreath product}). 
The hypotheses of Proposition~\ref{thm main} hold by the assumptions of the present corollary.  Therefore for every positive integer $n$ and for every non-trivial quotient $\bar A$ of $A^n$, the embedding problem (3) of Proposition~\ref{thm main} has no proper solution.
Hence
$\psi_1\nek \psi_n\in \Psi$ induce $\nu_1\nek \nu_n \in \calN$ which
are independent and proper. Therefore $\calN$ is a family of
independent proper solutions of cardinality $m$.
\end{proof}

\section{Semi-free subgroups}
\label{sec:proofofmainthm}

\subsection{Proof of \ref{thm subgroups}}

Let $F$ be semi-free of rank~$m$ and let $M \le F$.
\subsubsection{Case~\ref{opensubgroup}}
Assume that $M$ is open in $F$. We apply Corollary~\ref{cor_main}.
Given an open $D \normal F$, we take an open $L \normal F$ with
$L\leq M\cap D$. Then for $F_0 = M$ and $N = L$, there are no proper
solutions of the embedding problem appearing in Corollary~\ref{cor_main}, since $\phi$ is an isomorphism and
$\alpha$ is not. Therefore, $M$ is semi-free.

\subsubsection{Case \ref{finitelygenerated}}
Assume that $F/\Mhat$ is finitely generated, where $\Mhat =
\bigcap_{\sigma\in F} M^\sigma$ is the normal core of $M$ in $F$.

We apply Proposition~\ref{thm main}. Let $\calE_1(A) = (\mu\colon
M\to G_1, \alpha_1\colon A\rtimes G_1\to G_1)$ be a nontrivial FSEP
for $M$. Let $D$ be an open normal subgroup of $F$ with $M \cap D
\leq \Ker\mu$. Let $F_0 = MD$ and $N = \Mhat\cap D$. Then $F/N$ is
finitely generated (as an open subgroup of $F/\Mhat\times F/D$).
Thus, $F$ has only finitely many open subgroups containing $N$ of
index at most $r = (F:D)|A|^{2}$. Their intersection, $L$, is an
open normal subgroup of $F$ containing $N$ and contained in $D$.

Now, for $n=2$, the embedding problem \eqref{eq twisted wreath
product EP}, i.e.\
\[
  (\phigag\colon F/N\to F/L, \alphagag\colon\Agag\wr_{F_0/L}F/L \to F/L),
\]
for any nontrivial quotient $\Agag$ of $A^2$, has no proper
solution. Indeed, assume there exists a proper solution
$\psigag\colon F/N\to \Agag\wr_{F_0/L} F/L$ of \eqref{eq twisted
wreath product EP}. By Lemma \ref{lem large subgroup} there is a
subgroup $B$ of $H = \Agag\wr_{F_0/L}F/L$ of index $(H:B) =
(F:F_0)|\Agag| \le r$ that does not contain $\Ker\alphagag$. In
particular, $(H:B) > (H:B\,\Ker\alphagag) = (F/L : \alphagag(B))$.
Write $\psigag^{-1}(B)$ as $K/N$, for some $N\leq K \le F$. Then
$(F:K) = (F/N : K/N) = (H : B) \le r$, and hence $L \le K$. As
$\phigag = \alphagag\circ \psigag$, we have $ K/L = \phigag (K/N) =
\alphagag(\psigag(K/N)) = \alphagag(B). $ Therefore
$$
(H:B) = (F:K) = (F/L:K/L) = (F/L: \alphagag(B)) < (H:B),
$$
a contradiction.

Since $F$ is semi-free, there exists a family $\Psi$ of independent,
and in particular pairwise independent, proper solutions of the
nontrivial FSEP $\calE(A) = (\phi\colon F\to F/L, \alpha\colon
A\wr_{F_0/L} F/L \to F/L)$ such that $|\Psi| = m$. By
Proposition~\ref{thm main}(\ref{induce independent}) with $n=2$,
$\Psi$ induces a family $\calN$ of pairwise independent proper
solutions of $\calE_1$ and $|\calN|=|\Psi|=m$. By Corollary~\ref{cor
pairwise} we get that $M$ is semi-free of rank $m$.

\subsubsection{Cases~\ref{weissauer}, \ref{diamond}, and \ref{pronilpotent}}
The proof of Case \ref{diamond} is verbally identical with the proof
of the Diamond Theorem, \cite[Theorem 25.4.3]{FriedJarden2005},
provided that we replace \cite[Proposition 25.4.1]{FriedJarden2005}
by our Corollary~\ref{cor_main}.

Case \ref{weissauer} immediately follows from Case \ref{diamond}. So
does Case \ref{pronilpotent}: Since $(F:M) = (F/N:M/N)$ is divisible
by two primes and the Sylow subgroups are normal in $F/N$, there are
two (Sylow) normal subgroups $P_1,P_2$ of $F/N$ such that $P_1\cap
P_2 = 1$ and $P_1,P_2\not\subseteq M/N$. The preimages $M_1,M_2$ of
$P_1,P_2$ are normal in $F$ and satisfy $M_1\cap M_2 = N \leq M$,
but $M_1\not\leq M$ and $M_2\not\leq M$.

\subsubsection{Case~\ref{commutator}}
Assume that $M \normal F$ and $F/M$ is abelian. It follows that $M$
is also semi-free either by Cases \ref{finitelygenerated} and
\ref{diamond} or directly from Corollary~\ref{cor_main}. We show the
former. If $F/M$ is cyclic, then, by Case \ref{finitelygenerated},
$M$ is semi-free. Otherwise, there exists a pro-$p$ subgroup of rank
$2$ in $F/M$, say $H$. It factors as $H = C_1\times C_2$, where
$C_1,C_2$ are nontrivial cyclic pro-$p$ group. Then $C_1\cap C_2 =1$
and $C_1,C_2\normal F/M$ (since $F/M$ is abelian). The preimages
$M_1,M_2$ of $C_1,C_2$ are normal in $F$ and satisfy $M_1\cap M_2 =
M$, but $M_1\not\leq M$ and $M_2\not\leq M$.

\subsubsection{Cases~\ref{weight}, \ref{sparse}, and \ref{finite exponent}}
\label{sec:weight}
The proofs of these three cases are based on Case~\ref{opensubgroup}
and on more elementary arguments than the other cases.

Recall  that $\weight(F/M)=1$ if $M$ is open, and $\weight(F/M)$ is
the cardinality of the set of open subgroups of $F$ that contain $M$
if $(F:M) = \infty$ (\cite[Section 25.2]{FriedJarden2005}).

\begin{proof}[Proof of Case~\ref{weight}]
Let $\calE(M)= (\phi\colon M\to G, \alpha\colon H\to G)$ be a FSEP
for $M$ and let $M_0 = \Ker \phi$. There is an open $D\normal F$
such that $D\cap M \le M_0$. By Case~\ref{opensubgroup} we may
replace $F$ by its open subgroup $D M$ to assume that $D M = F$.
Then $d m \mapsto \phi(m)$, for $d \in D, \ m \in M$, extends $\phi$
to an epimorphism $\phi\colon F\to G$. Let $F_0$ be its kernel. It
contains $D$, hence $F_0 M = F$ and $F_0 \cap M = M_0$. Thus
$(M:M_0) = (F:F_0)$ and we have the FSEP $\calE(F)=(\phi\colon F\to
G, \alpha\colon H\to G)$.

Let $\Psi$ be a family of independent proper solutions of $\calE(F)$
of cardinality $m$. Each $\psi\in \Psi$ defines a solution $\psi' :=
\psi|_{M}$ of $\calE(M)$. Let $\Psi' = \{\psi' \st \psi\in \Psi\}$
and let $\calX\subseteq \Psi'$ be a maximal subset of independent
proper solutions (Zorn's Lemma). We claim that $\calX$ has
cardinality~$m$.

Assume differently, that is to say, assume $|\calX| < m$. Let $N =
\bigcap_{\psi'\in \calX} \Ker\psi'$ if $\calX \ne \emptyset$ and $N
= M_0$ if $\calX = \emptyset$. In both cases $N \le M_0$.

It suffices to find $\psi\in \Psi$ such that $N\Ker\psi = F_0$.
Indeed, then for every open subgroup $N_0$ of $M_0$ containing $N$
we have $(N_0:N_0\cap \Ker\psi)=(F_0:\Ker\psi)$,
\[
\xymatrix{
&& M \ar@{-}[d] \ar@{-}[r] & F \ar@{-}[d] \\
N \ar@{-}[d] \ar@{-}[r] & N_0 \ar@{-}[d] \ar@{-}[r] & M_0 \ar@{-}[d] \ar@{-}[r] & F_0 \ar@{-}[d] \\
N\cap \Ker\psi \ar@{-}[r] & N_0\cap \Ker\psi \ar@{-}[r] & M\cap \Ker\psi=\Ker\psi' \ar@{-}[r] & \Ker\psi \\
}
\]
i.e., $N_0$ and $\Ker\psi'$ are $M_0$-independent. In particular,
taking $N_0 = M_0$, we have $(M_0:\Ker\psi') = (M_0:M\cap
\Ker\psi)=(F_0:\Ker\psi)$, and hence $\psi'$ is surjective.
Furthermore, for any finite subset $\calX'$ of $\calX$, taking $N_0
= \bigcap_{\psi'\in \calX'} \Ker\psi'$ we get by
Proposition~\ref{proposition independent}(\ref{inductive
independent}) that $\calX'\cup\{\psi'\}$ is an independent set of
solutions. Therefore so is $\calX\cup\{\psi'\}$, which contradicts
the maximality of $\calX$.

To complete the proof, for each $\psi\in \Psi$ let $L_\psi =
N\Ker\psi$ and assume that $L_\psi \neq F_0$. Since $\{\Ker\psi \st
\psi\in \Psi\}$ is $F_0$-independent, the set $\{L_\psi \st \psi\in
\Psi\}$ is also independent by Proposition~\ref{proposition
independent}\eqref{intermediate}. Since $L_\psi \neq F_0$ for all
$\psi \in \Psi$, this implies in particular $L_{\psi_1}\neq
L_{\psi_2}$ for all distinct $\psi_1, \psi_2\in \Psi$. Hence
$\weight(F_0/N)\geq m$. But $\weight(F_0/M) < m$ by the hypothesis
of Case~\ref{weight} and the fact that $F_0$ is an open subgroup of
$F$. Moreover $\weight(M/N)<m$, by
\cite[Lemma~25.2.1(b)]{FriedJarden2005}.  Hence $\weight(F_0/N) < m$
by \cite[Lemma~25.2.1(d)]{FriedJarden2005}, a contradiction.
\end{proof}

\begin{definition}\label{def_sparse}
A closed subgroup $M$ of a profinite group $F$ of infinite index is
called \textbf{sparse} if for all $n\in\bbN$ there exists an open
subgroup $K$ of $F$ containing $M$ such that for every proper open
subgroup $L$ of $K$ containing $M$ we have $(K:L)\geq n$.
\end{definition}

The following lemma shows that this definition is equivalent to
\cite[Definition 2.1]{Bary-Soroker2006}:
\begin{lemma}
If $M$ is sparse in $F$, then for every $\ell,n\in \bbN$ there
exists $K$ as in Definition~\ref{def_sparse} of index at least
$\ell$ in $F$.
\end{lemma}

\begin{proof}
Let $\ell,n\in \bbN$. Choose an open subgroup $K_0$ of index
$\ell_0\geq \ell$ in $F$ such that $M\leq K_0$. By the definition
there exists $K_1$ with $M\leq K_1\leq F$ such that $(K_1:L) \geq
n\ell_0$ for all proper open subgroups $L$ of $K_1$ that contain
$M$. Then the assertion follows with $K= K_0\cap K_1$, since
$(K_1:K) \leq \ell_0$.
\end{proof}

\begin{proof}[Proof of Case~\ref{sparse}]
Let $M$ be a sparse subgroup of $F$. Let $\calE_0(A) = (\mu\colon M
\to G, \alpha\colon A\rtimes G\to G)$ be a nontrivial FSEP for $M$.

Choose an open normal subgroup $E_0$ of $F$ such that $E_0 \cap M
\leq \Ker \mu$ and let $F_0 = ME_0$. Since $M$ is sparse in $F_0$
\cite[Corollary 2.3]{Bary-Soroker2006}, there is an open subgroup
$K$ of $F_0$ containing $M$ such that $(K:L) > |A|^2|G|$ for each
proper open subgroup $L$ of $M$ that contains $M$. Extend $\mu$ to
an epimorphism $\phi\colon K\to G$ by $\phi(re) = \mu(r)$, $r\in M$,
$e\in E_0$. By Case~\ref{opensubgroup}, $K$ is semi-free of rank
$m$; hence it suffices to show that two independent proper solutions
$\psi_1,\psi_2$ of $\calE(A) = (\phi\colon K\to G, \alpha \colon A
\rtimes G \to G)$ induce two independent proper solutions
$\psi_1|_M,\psi_2|_M$ (Corollary~\ref{cor pairwise}).

By Lemma~\ref{lem twisted fiber}, $A^2\rtimes G$ is the fiber
product of $A\rtimes G\to G$ with itself. Thus $\psi_1,\psi_2$
induce a proper solution $\psi$ of $\calE(A^2)=(\phi\colon K\to G,
\alpha \colon A^2\rtimes G \to G)$ (Lemma~\ref{lem_fiber
independent}). Let $L = \Ker\psi$. Then $(K:ML) = (A^2 \rtimes G :
\psi(M)) \le |A|^2|G|$. Hence, by the choice of $K$, we get that $ML
= K$. Therefore, $\psi|_{M}$ is a proper solution of $\calE_0(A^2) =
(\phi\colon M\to G, \alpha \colon A^2\rtimes G \to G)$. But
$\psi|_{M}=\psi_1|_M \times \psi_2|_M$. Consequently,
$\psi_1|_M,\psi_2|_M$ are independent proper solutions of
$\calE_0(A)$, as claimed.
\end{proof}

The following corollary of Case~\ref{sparse} extends \cite[Lemma
2.4]{Bary-Soroker2006} to free groups of uncountable infinite rank.

\begin{corollary}\label{cor sparse implies finite
exponent} If $M$ is a sparse subgroup of a free profinite group $F$
of rank $m\geq 2$, then $M$ is a free profinite group of $\rank(M) =
\max\{\aleph_0, \rank(F)\}$.
\end{corollary}

\begin{proof}
The case where $\rank(F)\leq \aleph_0$ is proven in
\cite{Bary-Soroker2006}. Assume $m = \rank(F)$ is infinite. By
Theorem~\ref{thm QF and Proj imply Free}, $F$ is semi-free of
rank~$m$. By Case~\ref{sparse} of the \themainthm, $M$ is semi-free
of rank $m$.  Also, $M$ is projective, being a closed subgroup of a
free profinite group. Consequently $M$ is free of rank $m$
(Theorem~\ref{thm QF and Proj imply Free}).
\end{proof}

Case~\ref{finite exponent} is, in fact, a special case of
Case~\ref{sparse}:

\begin{lemma}\label{new sparse group}
Let $M$ be a closed subgroup of a profinite group $F$ of infinite
index. Assume $(F:M) = \prod_p p^{\alpha(p)}$ with all $\alpha(p)$
finite. Then $M$ is sparse in $F$.
\end{lemma}

\begin{proof}
For $n\in \bbN$ take $K$ to be an open subgroup of $F$ containing
$M$ such that $p^{\alpha(p)}\mid (F:K)$ for all $p\leq n$. Then for
each $M\leq L \lneqq K$ only primes $p > n$ can divide $(K:L)$.
Therefore, $(K:L) > n$.
\end{proof}

As a consequence of Corollary~\ref{cor sparse implies finite
exponent} and Lemma~\ref{new sparse group}, we get \cite[Proposition
5.1]{JardenLubotzky1992}:

\begin{corollary}
Let $M$ be a closed subgroup of a free profinite group $F$ of rank
$m\geq 2$. Assume $(F:M) = \prod_p p^{\alpha(p)}$ with all
$\alpha(p)$ finite.  If $(F:M)$ is infinite, then $M$ is free
profinite group of rank $\max\{\aleph_0,\rank(F)\}$.
\end{corollary}

\section{Quasi-freeness vs.\ semi-freeness}\label{sec:qfvssqf}
We now construct an example of a quasi-free group that is not
semi-free.

For a profinite group $C$ and an infinite set $X$ denote by $\freeprod_XC$
the free product of copies $\{C_x\}_{x \in X}$ of $C$ in the sense
of \cite{BinzNeukirchWenzel1971}. That is, $\freeprod_XC$ contains a copy
$C_x$ of $C$ for each $x \in X$; and every family of homomorphisms
$\psi_x \colon C_x \to A$ into a finite group $A$, such that
$\psi_x(C_x) = 1$ for all but finitely many $x \in X$, uniquely
extends to a homomorphism $\psi \colon \freeprod_XC \to A$.  As usual let
$\Fhat_\omega$ denote the free profinite group of countable rank.

\begin{proposition}
\label{prop QF not SQF}%
Let $X$ be a set of infinite cardinality $m$. Let $C = \prod_p
\bbZ/p\bbZ$ be the direct product of all prime cyclic groups. Let $F
= (\freeprod_XC) * \Fhat_\omega$. Then
\begin{enumerate}
\item
\label{F is QF} $F$ is quasi-free of rank~$m$,  and
\item
\label{F is not SQF} the FSEP
\begin{equation}\label{eq C_4 EP}
 (F\to 1, \bbZ/4\bbZ \to 1)
\end{equation}
has at most countably many independent proper solutions.
\end{enumerate}
In particular, for $m > \aleph_0$, $F$ is quasi-free but not
semi-free.
\end{proposition}

\begin{proof}
(\ref{F is QF}) The rank of $\freeprod_XC$ is $m$ and the rank of
$\Fhat_\omega$ is $\aleph_0 \le m$. Hence the rank of $F$ is $m$. In
particular, every FSEP for $F$ has at most $m$ proper solutions. Let
\begin{equation}
(\phi \colon F \to G, \alpha \colon H \to G) \label{eq}
\end{equation}
be a nontrivial FSEP. Let $\beta\colon G \to H$ be its splitting. We
need two auxiliary maps: Firstly, there exists a nontrivial
homomorphism $\pi \colon C \to \Ker\alpha$; namely, an epimorphism
of $C$ onto a subgroup of $\Ker\alpha$ of prime order. Secondly,
since $\Fhat_\omega$ is free of infinite rank, there exists an
epimorphism $\psi' \colon \Fhat_\omega \to
\alpha^{-1}(\phi(\Fhat_\omega))$ such that $\alpha \circ \psi'$ is
the restriction of $\phi$ to $\Fhat_\omega$. In particular,
$\psi'(\Fhat_\omega)$ contains $\Ker\alpha$. Since $\phi$ is
continuous, there is a $Y \cont X$ such that $X \hefresh Y$ is
finite and $\phi(C_y) = 1$ for every $y \in Y$.

For every $y \in Y$ define a homomorphism $\psi_y \colon F \to H$ in
the following manner: Its restriction to $C_y \isom C$ coincides
with $\pi$; if $y \ne x \in Y$, the restriction of $\psi_y$ to $C_x$
is trivial; if $x \in X \hefresh Y$, the restriction of $\psi_y$ to
$C_x$ is $\beta\circ \phi$; and, finally, the restriction of
$\psi_y$ to $\Fhat_\omega$ is $\psi'$. Thus $\alpha \circ \psi_y =
\phi$. As $\psi_y(F) \supseteq \psi'(\Fhat_\omega) \supseteq
\Ker\alpha$, the map $\psi_y$ is a proper solution of \eqref{eq}.

As $\psi_{y_1} \ne \psi_{y_2}$ for distinct $y_1, y _2 \in Y$, \
\eqref{eq} has at least $|Y| = m$ distinct proper solutions.

(\ref{F is not SQF}) Let $\Psi$ be an independent set of proper
solutions of \eqref{eq C_4 EP}. The map $\alpha\colon\bbZ/4\bbZ \to 1$
decomposes as $\alpha = \beta\gamma$, where  $\gamma \colon \bbZ/4\bbZ \to \bbZ/2\bbZ$ and $\beta
\colon \bbZ/2\bbZ \to 1$. If $\psi_1, \psi_2 \in \Psi$ are
independent, then $\gamma \circ \psi_1, \gamma \circ \psi_2$ are
independent proper solutions of $(\beta \colon \bbZ/2\bbZ \to 1,
\phi \colon F \to 1)$ (Proposition~\ref{proposition
independent}(d)). In particular, $\gamma \circ \psi_1 \ne \gamma
\circ \psi_2$. Thus $\{\gamma \circ \psi \st \psi \in \Psi\}$ has at
least the cardinality of $\Psi$.

On the other hand, $\bbZ/4\bbZ$ is a $2$-group and the $2$-Sylow
subgroup of $C$ is of order~$2$. Hence every $\psi \in \Psi$ maps
each $C_x \isom C$ into  $\Ker\gamma$, the unique subgroup of
$\bbZ/4\bbZ$ of order $2$, and hence $\gamma\circ \psi$ is trivial
on $C_x$. Therefore $\gamma\circ \psi$ is trivial on $\freeprod_XC$. It
follows that $\gamma\circ \psi$ is determined by its restriction to
$\Fhat_\omega$. But there are $\aleph_0$ (continuous) homomorphisms
$\Fhat_\omega \to \bbZ/4\bbZ$. Thus $|\Psi| \le \aleph_0$.
\end{proof}

\begin{remark}
One can modify the construction in the proposition to get an
absolute Galois group $F$ which is quasi-free but not semi-free.
E.g., let $F = \freeprod(\prod_{p \ne 2} \bbZ_p) * D * \Fhat_\omega$,
where $D$ is the free product of the constant sheaf of copies of
$\bbZ/2\bbZ$ over some profinite space of weight $m$. One can show
along the lines of the proof of Proposition~\ref{prop QF not SQF}
that $F$ is quasi-free but not semi-free. Moreover, $F$ is real
projective in the sense of \cite[p.~472]{HaranJarden1985} and hence
isomorphic to an absolute Galois group by
\cite[Theorem~10.4]{HaranJarden1985}. We leave out the details,
since the assertion is outside the scope of this work.
\end{remark}

\begin{remark}
In order to complete the picture we show that being semi-free is strictly weaker than being free. In fact, if $F$ is semi-free of infinite rank $m$ and $G$ is of rank $\leq m$, then $F\ast G$ is semi-free. This leads to many examples of semi-free but not free profinite groups; e.g., take $G$ to be finite and recall that a free group has no torsion.  Furthermore, we can construct a semi-free group of arbitrary cohomological dimension $d$, by taking $F$ free and $G$ of cohomological $d$.  If $d>1$ then the group is not free, or even projective, since its cohomological dimension is greater than one.  Another example is the absolute Galois group given in Theorem~\ref{thm laurent series} below, which is semi-free but is not projective (and hence not free) because its cohomological dimension is greater than one.  
\end{remark}

The condition $m > \aleph_0$ in the above proposition is essential:

\begin{remark}
If $\rank(F) = \aleph_0$, then $F$ is semi-free if and only if it is
quasi-free.

Indeed, assume $F$ is quasi-free. Then every FSEP is solvable. By
Lemma~\ref{lem:rankaleph0} $F$ is semi-free. The opposite direction
is immediate.

\end{remark}

We now show that Case~\ref{weight} of our \themainthm does not
carry over to quasi-free subgroups of quasi-free groups.

\begin{example}
\label{example qf no permanence} Let $X$ be a set of cardinality
$m>\aleph_0$ and let $F = (\freeprod_XC) * \Fhat_\omega$ be the group of
Proposition~\ref{prop QF not SQF}. Let $M$ be the kernel of the map
$F \to \Fhat_{\omega}$. Then $F$ is quasi-free of rank $m$,
$\weight(F/M) < m$, but $M$ is not quasi-free.

Indeed, by Proposition~\ref{prop QF not SQF}, $F$ is quasi-free of
rank $m$. We have
\[
\weight(F/M) = \rank(\Fhat_\omega) = \aleph_0
\]
since $F/M =
\Fhat_\omega$. It is easy to see that $M$ is generated by the
conjugates of $\freeprod_XC$ in $F$. Since $\freeprod_XC$ is generated by copies
of $C$ and $C = \prod_p \bbZ/p\bbZ$ is generated by elements of
prime order, also $M$ is generated by elements of prime order. Hence
$\bbZ/q^2\bbZ$ is not an image of $M$. In particular, $M$ is not
quasi-free.
\end{example}

\begin{remark}
It is interesting to ask which of the cases of the \themainthm holds for quasi-free groups. As we have seen, Case~\ref{weight} does not hold. In \cite{RibesZalesskiiStevenson2006} Case~\ref{opensubgroup} is proved. Case~\ref{commutator} is proved in \cite{Harbater2006} for $M=[F,F]$. Combining the  methods of this paper together with \cite{Harbater2006}, one can extend the result to any $M$ such that $F/M$ is abelian but not a pro-$p$ group. 
The proof of Case~\ref{sparse} (and hence of \eqref{finite exponent}) can be carried over to quasi-free groups. 
However, we do not know if the diamond theorem, i.e.\ Case~\ref{diamond}, which is the central result of this paper, holds for quasi-free groups. All other cases are open in the quasi-free case. 

In order to use our method, i.e.\ using wreath products, for quasi-free groups for $M$ of infinite index in $F$, one needs to come up with a new idea, as explained at the end of Section 1.
\end{remark}

\section{Fields with semi-free absolute Galois groups}
\def \Gal {{\rm Gal}}
\def \s {{\rm s}}
\def \card {{\rm card\,}}

The main result in \cite{HarbaterStevenson2005} (Theorem~5.1 there)
was that for any field $k$, the absolute Galois group of $K :=
k((x,t))$ is quasi-free. In fact more is true:

\begin{theorem}\label{thm laurent series}
 Let $k$ be a field.
Then the absolute Galois group of the field $K := k((x,t))$ is
semi-free of rank ${\rm card}\, K$.
\end{theorem}

The proof of this stronger result is essentially contained in the
proof of the original theorem in \cite{HarbaterStevenson2005}. We
explain below what additional observations need to be made to
complete the argument, and how these observations also yield
stronger forms of other results in \cite{HarbaterStevenson2005}.  See also \cite[Theorem~5.1]{HarbaterStevenson2ed2005} for more details. 

First we recall the strategy used to prove
\cite[Theorem~5.1]{HarbaterStevenson2005}. The proof of that theorem
relied on a related geometric assertion,
\cite[Proposition~5.3]{HarbaterStevenson2005}. That proposition
asserted that given a split short exact sequence $1 \to N \to \Gamma
\ {\buildrel f \over \to} \ G \to 1$ of finite groups with
non-trivial kernel, any $G$-Galois connected normal branched cover
$Y^* \to X^* = \Spec\,k[[x,t]]$ can be dominated by a
$\Gamma$-Galois connected normal branched cover $Z^* \to X^*$.
Moreover it said that this cover may be chosen such that $Z^* \to
Y^*$ satisfied a splitting condition (that $Z^* \to Y^*$ is totally split at the generic points of the ramification locus of $Y^* \to X^*$), and that the set of
isomorphism classes of such covers $Z^* \to X^*$ has cardinality
equal to $m := \card k((x,t))$.

The proof of \cite[Proposition~5.3]{HarbaterStevenson2005} relied on
\cite[Theorem~4.1]{HarbaterStevenson2005}, which was a more global
version of that assertion. Namely, it considered a smooth connected
curve $X$ over a field $\hat k := k((t))$, and then considered a
finite split embedding problem for the absolute Galois group of the
function field $K$ of $X$ (this field $K$ being a global analog of
the more local field $K$ considered in
\cite[Proposition~5.3]{HarbaterStevenson2005}). The conclusion was
similar: that any $G$-Galois branched cover $Y \to X$ of normal
curves can be dominated by a $\Gamma$-Galois branched cover $Z \to
X$; that this cover can be chosen with a splitting property; and that there are $m := \card
K$ distinct such choices of corresponding normal branched covers $Z
\to X$.  (The splitting property is that $Z \to Y$ is totally split over a given finite set $D \subset Y$ of closed points, and the decomposition groups of $Z \to X$ at the points of $Z$ over $\delta \in D$ are the conjugates of $\sigma(G_\delta)$, where $G_\delta$ is the decomposition group of $Y \to X$ at $\delta$ and where $\sigma$ is a section of $f$.)

Moreover, for the sake of
\cite[Proposition~5.3]{HarbaterStevenson2005}, more was shown in
\cite[Theorem~4.1]{HarbaterStevenson2005}, to enable passage from a
global solution to a more local solution.  Let $\bar X$ be a smooth
projective model for $X$ over $k[[t]]$; and with $Y,Z$ as above, let
$\bar Y, \bar Z$ be the corresponding normal branched covers. Let
$P$ be a closed point of $\bar X$ whose residue field is separable
over $k$, let $X^*$ be the spectrum of the complete local ring of
$\bar X$ at $P$, and suppose that the pullback $Y^* \to X^*$ of
$\bar Y \to \bar X$ is connected. Then among the pullbacks $Z^* \to
X^*$ of the above solutions $\bar Z \to \bar X$ there are $m$
distinct proper solutions of the corresponding local embedding
problem. This additional condition was applied in the case of the
$x$-line over $\hat k$ in order to obtain
\cite[Proposition~5.3]{HarbaterStevenson2005}.

More specifically, the relationship between the local assertion \cite[Proposition~5.3]{HarbaterStevenson2005} and the more global assertion \cite[Theorem~4.1]{HarbaterStevenson2005} is based on viewing $k((x,t))$ as the fraction field of the complete local ring of $\bar X := \bbP^1_{k[[t]]}$ at the point $x=t=0$.  In order to apply \cite[Theorem~4.1]{HarbaterStevenson2005} to the proof of \cite[Proposition~5.3]{HarbaterStevenson2005}, a change of variables can be made to reduce to the case in which the prime $(t)$ is unramified in $Y^* \to X^*$.  The reduction of this cover modulo $(t)$ is then induced from a branched cover of the projective $k$-line, by the Katz-Gabber theorem \cite[Theorem 1.4.1]{Katz1986}.  A patching argument then shows that this cover of $\bbP^1_k$ is in turn the closed fiber of a cover of $\bbP^1_{k[[t]]}$ that restricts to $Y^* \to X^*$.  This enables \cite[Theorem~4.1]{HarbaterStevenson2005} to be cited; and by the extra conditions in the paragraph above, the proper solutions to the embedding problem over the function field of $\bbP^1_{k[[t]]}$ yield distinct proper solutions to the embedding problem over $k((x,t))$.

Theorem~4.1 of \cite{HarbaterStevenson2005} was a variant on results
of Pop \cite[Main Theorem~A]{Pop1996} and of Haran and Jarden
\cite[Theorem~6.4]{HaranJarden1998}, showing that finite split
embedding problems over the function fields of curves over complete
discretely valued (or more generally large) fields have proper
regular solutions (and that some additional conditions can also be
satisfied, e.g.\ the existence of an unramified rational point).
Like those earlier results,
\cite[Theorem~4.1]{HarbaterStevenson2005} was proven using patching.
Generators were chosen for the kernel $N$ of the given finite split
embedding problem; and cyclic covers were constructed with groups
generated by each of those elements in turn. These were then patched
together to form a global solution; in doing so, a compatibility
condition (agreement on overlaps) had to be satisfied by the cyclic
covers on the ``patches''. Such a construction was carried out in
\cite[Proposition~3.5]{HarbaterStevenson2005}. But the construction
there assumed that branch points of $Z \to Y$ that correspond to
distinct generators of $Z$ had the property that their closures in
$\bar Y$ are disjoint. In order to apply this to the proof of
\cite[Theorem~4.1]{HarbaterStevenson2005} (where the branch points
all coalesce on the closed fiber at $P$, in order to preserve the
solutions over $X^*$), it was necessary to blow up the closed fiber
to separate the branch points.

We can now describe the proof of Theorem~\ref{thm laurent series}:

\begin{proof} As discussed above, this theorem is a strong form of~\cite[Theorem~5.1]{HarbaterStevenson2005}, and to prove this result
it suffices to prove a corresponding strong form
of~\cite[Proposition~5.3]{HarbaterStevenson2005}: that among the
covers $Z^* \to X^*$ whose existence is asserted in that
proposition, there is a subset having cardinality $m$, and which is
linearly disjoint as a set of covers of $Y^*$.  To prove this, we
need to see that in the situation of
\cite[Theorem~4.1]{HarbaterStevenson2005}, an additional property
holds: that there are $m$ choices of $Z \to X$ that are linearly
disjoint over $Y$, that properly solve the given global embedding
problem, and that induce proper solutions over $X^*$ that are
linearly disjoint over $Y^* = Y \times_X X^*$.

To show this stronger version of
\cite[Theorem~4.1]{HarbaterStevenson2005}, the key point is that the
branch points associated to the generators of $N$ can be chosen in
$m$ different (and even disjoint) ways. As shown in the original
proof, given any choices of these points on $X$ (which correspond to
curves on $\bar X$ that are finite over $k[[x]]$), any other choice of
points that is congruent to the original choice modulo a
sufficiently high power of $t$ will also work. (Indeed, this is how
it was shown that there are $m$ distinct solutions, both over $X$
and over $X^*$.)  What needs to be shown here is that by varying the
branch points we can obtain $m$ solutions that are linearly disjoint
over $Y$. Since Galois branched covers with no common subcover are
linearly disjoint, it suffices to show that the set of $m$ solutions
$Z \to X$, such that the covers $Z \to Y$ have pairwise disjoint
branch loci, can be chosen such that each $Z \to Y$ has no
non-trivial \'etale subcover $W \to Y$.

In the above situation, if $Z \to Y$ has a non-trivial \'etale
subcover $W \to Y$, then the Galois group $\Gal(Z/W)$, which is a
subgroup of $N = \Gal(Z/Y)$, must contain all the inertia groups of
$Z \to Y$. But this is ruled out by the explicit construction in the
proof of \cite[Proposition~3.5]{HarbaterStevenson2005}. Namely, that
result asserts that the closed fiber $\bar Z \to \bar Y$ of $Z \to
Y$ is an $N$-Galois mock cover; i.e., each irreducible component of
$\bar Z$ maps isomorphically onto $\bar Y$, with the irreducible
components being indexed by the cosets of $N$ in $\Gamma$. The
construction in the proof there shows that for each generator $n$ of
$N$, there is a closed point $Q_n \in \bar Z$ lying in the
ramification locus of $\bar Z \to \bar Y$, such that $n$ generates
the inertia group of $\bar Z \to \bar Y$ at $Q_n$ and also the
inertia groups at the generic points of the ramification components
passing through $Q_n$. Since the elements $n$ together generate $N$,
this shows that the $N$-Galois cover $Z \to Y$ has no non-trivial
\'etale subcovers, as desired.

Thus the above strong form of
\cite[Theorem~4.1]{HarbaterStevenson2005} indeed holds. Hence so
does the strong form of
\cite[Proposition~5.3]{HarbaterStevenson2005}; and thus also
Theorem~\ref{thm laurent series} above, the strong form of
\cite[Theorem~5.1]{HarbaterStevenson2005}.
\end{proof}

Another key result of \cite{HarbaterStevenson2005},
viz.~Corollary~4.4 there, asserted that if $K$ is the function field
of a smooth projective curve over a very large field $k$, then the
absolute Galois group of $K$ is quasi-free. This can also be
strengthened, as follows:

\begin{theorem}\label{thm curves over large fields}
If $K$ is the function field of a smooth projective curve $X_0$ over
a large field $k$, then the absolute Galois group of $K$ is
semi-free.
\end{theorem}

This result has been independently proved by Jarden \cite{Jarden}.

\begin{proof}
By a recent result of Pop (see
\cite[Proposition~3.3]{Harbater2006}), every large field is very
large. So the assumption on $k$ in
\cite[Corollary~4.4]{HarbaterStevenson2005} can be (a priori)
weakened from very large to large. Concerning the strengthening of
the conclusion, this can be done in a similar way to what was done
above for Theorem~\ref{thm laurent series}. Namely,
\cite[Corollary~4.4]{HarbaterStevenson2005} followed from
\cite[Theorem~4.3]{HarbaterStevenson2005}, which was a variant of
\cite[Theorem~4.1]{HarbaterStevenson2005} in which the field $\hat k
= k((t))$ was replaced by a more general large field $F$. As in the
case of Theorem~\ref{thm laurent series}, to prove \ref{thm curves
over large fields} it suffices to show that the proper solutions
$Z_0 \to X_0$ in \cite[Theorem~4.3]{HarbaterStevenson2005} can be
chosen so as to be linearly disjoint over $Y_0$; and for this it
suffices to show that they can be chosen so that each $Z_0 \to Y_0$
has no non-trivial \'etale subcovers.

Theorem~4.3 of \cite{HarbaterStevenson2005} was proven using
\cite[Theorem~4.1]{HarbaterStevenson2005}, by taking $k=F$;
obtaining a proper solution for the function field of the induced
curve $\bar X := X_0 \times_F R$ over $R = k[[t]]$; descending from
$R$ to a $k$-algebra $A$ of finite type, corresponding to a
$k$-variety $V$; considering the descended $\Gamma$-Galois cover
$Z_A \to X_A$ as a family of $\Gamma$-Galois covers of $X_0$
parametrized by $V$; and then specializing to $k$-points of $V$
(thereby obtaining solutions over $X_0$) using that $k$ is (very)
large. To prove the desired strong form of
\cite[Theorem~4.3]{HarbaterStevenson2005}, observe that in the
context of the above use of
\cite[Theorem~4.1]{HarbaterStevenson2005}, the branch points (which
can be varied arbitrarily modulo some sufficiently high power of
$t$) can be chosen so as not to be constant; i.e.\ not of the form
$P' \times_k \hat k$ with $P'$ a point of $X_0$. As a result, the
the varying branch locus of the family of $\Gamma$-Galois covers of
$X_0$ parametrized by $V$ is base-point free. So as in the proof of
the strong form of \cite[Theorem~4.1]{HarbaterStevenson2005}, the
specialized covers can be chosen to have no non-trivial \'etale
subcovers; and hence they are linearly disjoint. This shows that
\cite[Theorem~4.3]{HarbaterStevenson2005} can be strengthened as
claimed to include the desired linear disjointness assertion; and
hence Theorem~\ref{thm curves over large fields}, the strong form of
\cite[Corollary~4.4]{HarbaterStevenson2005}, also holds.
\end{proof}

\section{Fields with free absolute Galois groups}\label{sec free
absolute} We present two families of fields having free absolute
Galois groups. For each we use Theorem~\ref{thm QF and Proj imply
Free} to reduce the proof of freeness to proving that the group is
semi-free and projective.

The semi-freeness follows from the Diamond Theorem (\ref{thm
subgroups}, Case~\ref{diamond}) together with the semi-freeness of
the absolute Galois group of the base field, which was established
in the previous section. The projectivity is achieved by different
means (here we just quote it).

\subsection{Fields containing the maximal abelian extension of $k((x,t))$}
We follow \cite{Harbater2006} to find fields with free absolute
Galois group. Let us start with a general fact and then give some
concrete examples.

\begin{corollary} \label{cor:free}
Let $K = k((x,t))$, where $k$ is separably closed and let $L$ be a
separable extension of $K$. If $L$ contains the maximal abelian
extension of $K$, and its absolute Galois group $\gal(L)$ satisfies
one of the cases of the \ref{thm subgroups} as a subgroup of
$\gal(K)$, then $\gal(L)$ is a free profinite group.
\end{corollary}

\begin{proof}
The group $\gal(K)$ is semi-free of rank $m$ by Theorem~\ref{thm
laurent series}. Hence so is $\gal(L)$. Also, $\gal(L)$ is
projective \cite[Theorem 4.4]{Harbater2006} (see also
\cite{Colliot-TheleneOjangurenParimala2002}). Thus, Theorem~\ref{thm
QF and Proj imply Free} yields that $\gal(L)$ is free.
\end{proof}

\begin{example}
Let $K = k ((x,y))$, where $k$ is separably closed. Let $E$ be a
Galois extension of $K$ not containing the maximal abelian extension
$K^{\ab}$ of $K$. Let $L$ be any subextension of $E K^{\ab} /
K^{ab}$. We claim that $\gal(L)$ is free of rank equal to the
cardinality of $L$.

To see this, first note that $\gal(K)$ is semi-free
(Theorem~\ref{thm laurent series}). If $L = K^{\ab}$, then by
\cite[Theorem~4.6(b)]{Harbater2006} it follows that $\gal(L)$ is
free. (Equivalently, this follows from \ref{thm subgroups}
Case~\ref{commutator} together with Corollary~\ref{cor:free}.)

Now consider the case $L\neq K^{\ab}$. Since $K^{\ab}\not\subseteq
E$ and $K^{\ab}\subseteq L$, it follows that $L\not \subseteq E$.
Furthermore, $E/K$ and $K^{\ab}/K$ are Galois. Hence by the Galois
correspondence, $M=\gal(L)$ satisfies Case~\ref{diamond} of the Main
Theorem with $F=\gal(K)$, $M_1 = \gal(E)$, and $M_2 =
\gal(K^{\ab})$. By Corollary~\ref{cor:free}, $\gal(L)$ is free.
\[
\xymatrix{%
\gal(K^{\ab})\ar@{-}[d]|{\gal(L)}
    &\gal(K)\ar@{-}[l]\ar@{-}[d]\\
\gal(E)\cap \gal(K^{\ab})
    &\gal(E)\ar@{-}[l]
}%
\]
\end{example}

\subsection{Jarden's example -- extension of roots}
This example is adapted from \cite{Jarden}. Let $k$ be a PAC field
of characteristic $p \ge 0$ and $K = k(x)$. Let $\mathcal F\subseteq
k[x]\subseteq K$ be the set of all monic irreducible polynomials.
For each $f \in \mathcal F$ choose a set of compatible roots
\[
\big\{f^{\frac1n} \ \big|\  p\nmid n\big\} \subseteq K_s.
\]
(Here compatible means that $(f^{\frac{1}{nn'}})^n =
f^{\frac{1}{n'}}$ for all $n,n'$ prime to $p$.) Let
\[
L = K\big(f^{\frac{1}{n}} \ \big|\  f\in \mathcal F\ \mbox{and}\ p\nmid n\big).
\]
Note that $L/K$ is Galois if and only if $K$ contains all roots of
unity. Thus in general $L/K$ is not Galois. In what follows we show
that $\gal(L)$ is free of rank equal to the cardinality of $L$.

\begin{fact} \label{gal projectivity}
$\gal(L)$ is projective.
\end{fact}

This fact follows from theorems of Efrat and Pop (see Theorems~10.4.9
and~11.6.4 in \cite{Jarden}).

\begin{lemma} \label{diamond roots of polynomials}
There exist Galois extensions $L_1,L_2$ of $K$ such that $L\subseteq
L_1L_2$, but $L \not\subseteq L_i$, $i=1,2$.
\end{lemma}

\begin{proof}
Let $L_0$ denote the extension of $K$ generated by all roots of
unity. Let
\[
L_1 = L_0\big(x^{\frac1n}\ \big|\  p\nmid n)\ \mbox{and}\ %
L_2 = L_0\big(f^{\frac1n}\ \big| \ f\in \mathcal F\smallsetminus \{x\}\
\mbox{and}\ p\nmid n).
\]
Clearly $L_1,L_2$ are Galois extensions of $K$. It is obvious that
$L\subseteq L_1L_2$. Choose an integer $m>1$ that is not
divisible by $p$. Since $(x+1)^{\frac{1}{m}}\not \in L_1$ we get
that $L\not \subseteq L_1$; and similarly $x^{\frac{1}{m}}\not\in
L_2$ implies that $L\not \subseteq L_2$.
\end{proof}

\begin{theorem} \label{freeness for root field}
$\gal(L)$ is free of rank equal to the cardinality of $L$.
\end{theorem}

\begin{proof}
By Theorem~\ref{thm QF and Proj imply Free} it suffices to show that
$\gal(L)$ is both projective and semi-free of rank equal to the
cardinality of $L$. We already mentioned that $\gal(L)$ is
projective (Fact~\ref{gal projectivity}).

Theorem~\ref{thm curves over large fields} implies that $\gal(K)$ is
semi-free of rank $m := |K| = |L|$. (Recall that $k$ is PAC, and in
particular large.) Taking absolute Galois groups of the fields
$L_1,L_2$ in the above lemma establishes the condition of
Case~\ref{diamond} of the \themainthm, thus $\gal(L)$ is semi-free
of rank $m$.
\end{proof}

In fact, even more is true.  Namely, we have learned from Pop that
the proof of his theorem (referred to above) applies more broadly.
In particular, it applies in the case that $k = F((t))$ for some
separably closed field $F$ (using that this field $k$, like a PAC
field, has projective absolute Galois group and ``satisfies a
universal local-global principle'').  Following the same
construction as above, we again deduce that the resulting field $L$
has free absolute Galois group of rank $|L|$.  Note that by
Corollary~25.4.8 of \cite{FriedJarden2005}, this also implies that
the absolute Galois group of $F((t))(x)^\ab$ is free for $F$
separably closed.

Moreover, if $k'$ is the field obtained from $k$ by adjoining a set
of compatible $n^{\rm th}$ roots to all the non-zero elements of
$k$, then Pop's argument also shows that $L':=Lk'$ has projective
absolute Galois group in the case that $k$ is a local field such as
$\bbF_p((t))$ or $\bbQ_p$.  (Here the adjunction of additional roots
is to deal with the fact that $\Gal(k)$ is no longer projective.)
Since Lemma~\ref{diamond roots of polynomials} then holds with $L$
replaced by $L'$ (and with $L_i$ in the proof replaced by its
compositum with $k'$), the above proof of Theorem~\ref{freeness for
root field} then shows that $\gal(L')$ is a free profinite group.

\medskip\noindent\textbf{Acknowledgment}
We thank Moshe Jarden for the suggestion to
consider Case~\ref{weight} of the \themainthm.


%
%
\def\cprime{$'$}
\providecommand{\bysame}{\leavevmode\hbox
to3em{\hrulefill}\thinspace}
\providecommand{\MR}{\relax\ifhmode\unskip\space\fi MR }
\providecommand{\MRhref}[2]{%
  \href{http://www.ams.org/mathscinet-getitem?mr=#1}{#2}
} \providecommand{\href}[2]{#2}

\end{document}